\documentclass{article}
\usepackage{graphicx} % Required for inserting images

\usepackage{amssymb, amsmath, amsthm}
\usepackage{amsfonts, bm}

\usepackage{float}
\usepackage[margin=2.6cm]{geometry}

\usepackage[hyphens]{url}
\usepackage[colorlinks=true,citecolor=black,linkcolor=black,urlcolor=blue]{hyperref}
\usepackage[hyphenbreaks]{breakurl}

\usepackage{xcolor}
\usepackage{enumitem}

\theoremstyle{plain}
\newtheorem{Th}{Theorem}
\newtheorem{Lemma}[Th]{Lemma}
\newtheorem{Cor}[Th]{Corollary}
\newtheorem{Prop}[Th]{Proposition}

\theoremstyle{definition}
\newtheorem{Def}[Th]{Definition}

\newtheorem{?}[Th]{Problem}

\newcommand{\doi}[1]{\url{https://doi.org/#1}}

\newcommand{\LW}{\mathrm{LW}}

\usepackage{authblk}

\title{Some bounds on the spectral radius of connected threshold graphs}

\author[1,2]{P\'eter Csikv\'ari}
\author[3,4,5]{Ivan Damnjanovi\'c}
\author[6]{Dragan Stevanovi\'c}
\author[7,8]{Stephan Wagner}

\affil[1]{HUN-REN Alfréd Rényi Institute of Mathematics, Budapest, Hungary}
\affil[2]{Eötvös Loránd University, Budapest, Hungary}
\affil[3]{Faculty of Mathematics, Natural Sciences and Information Technologies, University of Primorska, Koper, Slovenia}
\affil[4]{Faculty of Electronic Engineering, University of Niš, Niš, Serbia}
\affil[5]{Diffine LLC, San Diego, California, USA}
\affil[6]{College of Integrative Studies, Abdullah Al-Salem University, Khaldiya, Kuwait}
\affil[7]{Institute of Discrete Mathematics, TU Graz, Austria}
\affil[8]{Department of Mathematics, Uppsala University, Sweden}

\date{}

\begin{document}

\maketitle

\begin{abstract}
The spectral radius of a graph is the spectral radius of its adjacency matrix. A threshold graph is a simple graph whose vertices can be ordered as $v_1, v_2, \ldots, v_n$, so that for each $2 \le i \le n$, vertex $v_i$ is either adjacent or nonadjacent simultaneously to all of $v_1, v_2, \ldots, v_{i-1}$. Brualdi and Hoffman initially posed and then partially solved the extremal problem of finding the simple graphs with a given number of edges that have the maximum spectral radius. This problem was subsequently completely resolved by Rowlinson. Here, we deal with the similar problem of maximizing the spectral radius over the set of connected simple graphs with a given number of vertices and edges. As shown by Brualdi and Solheid, each such extremal graph is necessarily a threshold graph. We investigate the spectral radii of threshold graphs by relying on computations involving lazy walks. Furthermore, we obtain three lower bounds and one upper bound on the spectral radius of a given connected threshold graph.
\end{abstract}

\section{Introduction}

Throughout the paper, we will take all graphs to be undirected, finite and simple. The \emph{adjacency matrix} $A(G)$ of a graph $G$ is the binary $\mathbb{R}^{V(G) \times V(G)}$ matrix with zero diagonal such that $A(G)_{uv} = 1$ if and only if $u \sim v$ in $G$. The \emph{spectrum} $\sigma(G)$ of a graph $G$ is the spectrum of its adjacency matrix $A(G)$, i.e., the multiset of eigenvalues of $A(G)$. The \emph{spectral radius} $\rho(G)$ of a graph $G$ is defined as $\rho(G) = \max_{\lambda \in \sigma(G)} |\lambda|$. The Perron--Frobenius theory (see, e.g., \cite[Chapter 8]{godroy2001}) guarantees that $\rho(G)$ is necessarily an eigenvalue of $A(G)$. Many graph theorists have found interest in investigating various properties concerning the spectral radii of graphs (see, e.g., \cite{stevanovic2018} and the references therein).

A \emph{threshold graph} is a graph whose vertices can be ordered as $v_1, v_2, \ldots, v_n$, so that for each $2 \le i \le n$, vertex $v_i$ is either adjacent or nonadjacent simultaneously to all of the vertices $v_1, v_2, \ldots, v_{i - 1}$. These graphs were introduced by Chv\'atal and Hammer \cite{chvaham1977} in 1977 while constructing an algorithm connected to integer programming. Given a threshold graph $G$ of order $n \in \mathbb{N}$, it is not difficult to observe that if we arrange its vertices in descending order with respect to their degree, then the obtained adjacency matrix $A(G) \in \mathbb{R}^{n \times n}$ is of \emph{stepwise form}, i.e., for any $1 \le i < j \le n$ such that $A(G)_{ij} = 1$, we have $A(G)_{k \ell} = 1$ for any $1 \le k < \ell \le n$ such that $k \le i$ and $\ell \le j$. Furthermore, any graph must be a threshold graph whenever its vertices can be arranged in such a way that the obtained adjacency matrix is of stepwise form. For more on threshold graphs, please refer to \cite{mahpel1995} and the references therein.

In 1976, Brualdi and Hoffman \cite[p.\ 438]{berfourversot1978} initially posed the extremal problem of maximizing the spectral radius over the set of graphs with a prescribed number of edges. Later on, Brualdi and Hoffman \cite{bruhoff1985} proved in 1985 that each such extremal graph is necessarily a threshold graph and proceeded to show that for any $k \in \mathbb{N}$, the graphs attaining the maximum spectral radius among those with $\binom{k}{2}$ edges are precisely the graphs of the form $K_k + rK_1$ for some $r \in \mathbb{N}_0$. Further work on this matter was done by Friedland \cite{friedland1985, friedland1988} and Stanley \cite{stanley1987}. In a subsequent paper from 1988, Rowlinson \cite{rowlinson1988} finalized the solution to the given problem by showing that for each $k, t \in \mathbb{N}$ such that $t < k$, the extremal graphs among those with $\binom{k}{2}+t$ edges are the graphs whose only nontrivial component arises from $K_k$ by adding a new vertex and then joining it to $t$ vertices of $K_k$.

Here, we investigate the similar, but unresolved, problem of characterizing the \textbf{connected} graphs with a prescribed number of vertices and edges that attain the maximum spectral radius. Brualdi and Solheid~\cite{brusol1986} proved in 1986 that each such extremal graph must be a threshold graph. We also note that Simić, Li Marzi and Belardo \cite{silimabe2004} obtained the same result in 2004, albeit with a different strategy. It is natural to describe a threshold graph through a binary \emph{generating sequence} $a_1 a_2 \cdots a_n$, where $a_i = 1$ (resp.\ $a_i = 0$) indicates that vertex $v_i$ is adjacent (resp.\ nonadjacent) to all of $v_1, v_2, \ldots, v_{i - 1}$. Now, for any $k \in \mathbb{N}$, we will use $G_{p_1, p_2, \ldots, p_k}$ to denote the threshold graph with the generating sequence
\[
    \underbrace{11\cdots11}_{p_1\mbox{\scriptsize\ times}} \, \underbrace{00\cdots00}_{p_2\mbox{\scriptsize\ times}} \,
    \underbrace{11\cdots11}_{p_3\mbox{\scriptsize\ times}} \, \cdots \,
    \underbrace{11\cdots11}_{p_k\mbox{\scriptsize\ times}}
\]
if $k$ is odd, and
\[
    \underbrace{00\cdots00}_{p_1\mbox{\scriptsize\ times}} \, \underbrace{11\cdots11}_{p_2\mbox{\scriptsize\ times}} \,
    \underbrace{00\cdots00}_{p_3\mbox{\scriptsize\ times}} \, \cdots \,
    \underbrace{11\cdots11}_{p_k\mbox{\scriptsize\ times}}
\]
if $k$ is even. With this in mind, some of the current partial results on the spectral radius maximization for connected graphs of a given order and size can be succinctly stated through Table \ref{known_results}.

\begin{table}[h!]
\centering
\begin{tabular}{|l|*{3}{l|}}\hline
Size with respect to the order & Order scope & Extremal graph & Source(s) \\\hline
$m = n - 1$ & any feasible $n$ & $G_{n - 1, 1}$ & \cite{brusol1986, cosi1957, lope1973}
\\\hline
$m = n$ & any feasible $n$ & $G_{2, n - 3, 1}$ & \cite{brusol1986}
\\\hline
$m = n + 1$ & any feasible $n$ & $G_{2, 1, n - 4, 1}$ & \cite{brusol1986}
\\\hline
$m = n + 2$ & any feasible $n$ & $G_{3, n - 4, 1}$ & \cite{brusol1986}
\\\hline
$m = n + t$, where $t \ge 3$ & sufficiently large $n$ & $G_{t + 1, 1, n - 3 - t, 1}$ & \cite{brusol1986, cvetrow1988}
\\\hline
$m = n + \binom{k}{2} - 1$, where $k \ge 4$ & any feasible $n$ & $G_{k, n - 1 - k, 1}$ and/or $G_{\binom{k}{2}, 1, n - 2 - \binom{k}{2}, 1}$ & \cite{bell1991}
\\\hline
\begin{tabular}{@{}l@{}} $m = n + \binom{k}{2} - 2$, \\ so that $2n \le m < \binom{n}{2} - 1$ \end{tabular} & any feasible $n$ & $G_{2, k - 1, n - 2 - k, 1}$ & \cite{oleroydri2002}
\\\hline
\end{tabular}
\caption{Partial results for the spectral radius maximization problem over the set of connected graphs of order $n \in \mathbb{N}$ and size $m \in \mathbb{N}_0$.}
\label{known_results}
\end{table}

As shown in Table \ref{known_results}, the remaining case to be fully settled is when $m = n - 1 + \binom{k}{2} + t$ for some $k, t \in \mathbb{N}$ such that $3 \le k \le n - 2$ and $t < k$. Here, it is reasonable to expect that the solution is either $G_{k - t, 1, t, n - 2 - k, 1}$ or $G_{\binom{k}{2} + t, 1, n - 2 - \binom{k}{2} - t, 1}$ (or both). Our aim is to build upon the existing results by deriving several bounds on the spectral radius of a given connected threshold graph.

In Section \ref{sc_prel} we outline a known method for expressing the spectral radius of a graph via lazy walks and then introduce the auxiliary notation to be used later on. Afterwards, in Section~\ref{sc_recurrence} we derive a recurrence relation for computing certain lazy walks in a connected threshold graph. This formula is subsequently applied in Sections~\ref{sc_lower_bound} and \ref{sc_upper_bound} to obtain a lower and an upper bound, respectively, on the spectral radius of the given graph. Finally, in Section~\ref{sc_extra_bounds} we give two additional lower bounds whose derivation does not rely on lazy walks and then use a similar strategy to provide an alternative proof of the upper bound previously obtained in Section \ref{sc_upper_bound}.

\section{Preliminaries}\label{sc_prel}

A \emph{lazy walk} of length $k \in \mathbb{N}_0$ in a given graph $G$ is any sequence $v_0 v_1 v_2 \cdots v_k$ of vertices from $V(G)$ such that for each $1 \le i \le k$, we have $v_i = v_{i - 1}$ or $v_i \sim v_{i - 1}$ in $G$. Also, an empty sequence is considered to be an \emph{empty lazy walk}. It is clear (see, e.g., \cite[Chapter 1]{brouhae2012}) that for any $k \in \mathbb{N}_0$, the entry $\left( A(G) + I \right)^k_{uv}$ contains the number of lazy walks of length $k$ in $G$ whose starting and ending vertex are $u$ and $v$, respectively. With this in mind, we give the next folklore lemma.

\begin{Lemma}\label{lazy_lemma}
    Let $G$ be a connected graph and let $U_1$ and $U_2$ be any two nonempty subsets of $V(G)$. Furthermore, for each $k \in \mathbb{N}_0$, let $\Omega_k$ denote the number of lazy walks of length $k$ in $G$ whose starting and ending vertex belong to $U_1$ and $U_2$, respectively. Then we have
    \[
        \lim_{k \to +\infty} \sqrt[k]{\Omega_k} = 1 + \rho(G).
    \]
\end{Lemma}
\begin{proof}
    Let $n = |V(G)|$. The spectral decomposition (see, e.g., \cite[Chapter 8]{godroy2001}) of the real symmetric matrix $A(G)$ yields
    \[
        (A(G) + I)^k = \sum_{i = 1}^n (\lambda_i + 1)^k x_i^\intercal x_i \qquad (k \in \mathbb{N}_0),
    \]
    where $\lambda_1 \ge \lambda_2 \ge \cdots \ge \lambda_n$ are the eigenvalues of $A(G)$ with the corresponding eigenvectors $x_1, x_2, \ldots, x_n$. Thus, if $w_1$ and $w_2$ denote the $\{0, 1\}^{V(G)}$ indicator vectors of the sets $U_1$ and $U_2$, respectively, we get
    \[
        \Omega_k = w_1^\intercal \, (A(G) + I)^k \, w_2 = \sum_{i = 1}^n \left( w_1^\intercal x_i^\intercal x_i w_2 \right) (\lambda_i + 1)^k .
    \]
    Observe that $\lambda_1 + 1 = 1 + \rho(G)$, while the Perron--Frobenius theory tells us that $|\lambda_i + 1| < 1 + \rho(G)$ for every $2 \le i \le n$. Also, we have $w_1^\intercal x_1^\intercal x_1 w_2 > 0$ since the Perron--Frobenius vector is always positive for a connected graph. With all of this in mind, we obtain
    \begin{align*}
        \pushQED{\qed}
        \lim_{k \to +\infty} \sqrt[k]{\Omega_k} &= \lim_{k \to +\infty} \sqrt[k]{w_1^\intercal x_1^\intercal x_1 w_2 (1 + \rho(G))^k + \sum_{i = 2}^n \left( w_1^\intercal x_i^\intercal x_i w_2 \right) (\lambda_i + 1)^k}\\
        &= \lim_{k \to +\infty} (1 + \rho(G)) \sqrt[k]{w_1^\intercal x_1^\intercal x_1 w_2 + \sum_{i = 2}^n \left( w_1^\intercal x_i^\intercal x_i w_2 \right) \left( \frac{\lambda_i + 1}{1 + \rho(G)}\right)^k} = 1 + \rho(G) .\qedhere
    \end{align*}
\end{proof}

In the remainder of the paper, we will let $G$ be a connected threshold graph of order $n \ge 4$ and size $n - 1 < m < \binom{n}{2}$ with the generating sequence $a_1 a_2 \cdots a_n$ and the corresponding vertices $v_1, v_2, \ldots, v_n$. For the sake of brevity, we will write $\rho$ instead of $\rho(G)$. We will also say that $v_i$ is a \emph{type~1} (resp.\ \emph{type~0}) vertex if $a_i = 1$ (resp.\ $a_i = 0$). We trivially observe that $a_n = 1$, i.e., that $v_n$ must be a type~1 vertex. Without loss of generality, we may also take $v_1$ to be a type~1 vertex. Finally, for a given lazy walk $v_{i_0} v_{i_1} v_{i_2} \cdots v_{i_k}$, we define its \emph{signature} as the binary sequence $a_{i_0} a_{i_1} a_{i_2} \cdots a_{i_k}$.

Now, let $c \ge 3$ and $z \in \mathbb{N}$ denote the number of type 1 and type 0 vertices in $G$, respectively, so that $c + z = n$. We introduce the \emph{backwards zero position} (BZP) sequence as the nonincreasing tuple $(b_1, b_2, \ldots, b_z) \in \mathbb{N}^z$ where $b_i$ denotes the number of type 1 vertices that appear after the $i$-th type 0 vertex in the generating sequence. For example, the BZP sequence of $G_{k - t, 1, t, n - 2 - k, 1}$ is
\[
    (t+1, \underbrace{1, 1, \ldots, 1}_{n - 2 - k \mbox{\scriptsize\ times}}),
\]
while the BZP sequence of $G_{\binom{k}{2} + t, 1, n - 2 - \binom{k}{2} - t, 1}$ is
\[
    (\underbrace{2, 2, \ldots, 2}_{\binom{k}{2} + t \mbox{\scriptsize\ times}}, \underbrace{1, 1, \ldots, 1}_{n - 2 - \binom{k}{2} - t \mbox{\scriptsize\ times}}) .
\]
The generating sequence, and thus the threshold graph itself, can always be reconstructed from the given BZP sequence and the number of vertices. Indeed, the graph $G$ is obtained from the numbers $c, b_1, b_2, \ldots, b_z$ by starting from $K_c$ and then adding $z$ additional vertices so that the $i$-th vertex is connected to the first $b_i$ vertices of $K_c$. Also, note that $\binom{c}{2} + \sum_{i = 1}^z b_i = m$, as well as $1 \le b_i \le c - 1$ for any $1 \le i \le z$.

\section{Counting lazy walks}\label{sc_recurrence}

For each $k \in \mathbb{N}$, let $\LW_k$ denote the number of lazy walks of length $k - 1$ in $G$ whose starting and ending vertex are both a type 1 vertex. For convenience, we will also take that $\LW_0 = 1$. In the present section, we will derive a recurrence formula for computing the $\left(\LW_k\right)_{k \in \mathbb{N}_0}$ sequence. To begin, for some $p \in \mathbb{N}_0$, let $S$ be a lazy walk signature that contains $p + 1$ ones interleaved with precisely $p$ blocks of arbitrarily many zeros. In other words, let $S$ be of the form
\begin{equation}\label{f_form}
    1\underbrace{0\cdots 0}_{\mbox{\scriptsize\ block } 1}1\underbrace{0\cdots 0}_{\mbox{\scriptsize\ block } 2}1\underbrace{0\cdots 0}_{\mbox{\scriptsize\ block } 3}1 \cdots 1\underbrace{0\cdots 0}_{\mbox{\scriptsize\ block } p}1 .
\end{equation}
We claim that the number of lazy walks with the given signature $S$ is the same for any such $S$. Indeed, it does not matter how many zeros there are in each block, since the nonadjacency of type 0 vertices implies that in every lazy walk all the zeros from the same signature block correspond to the same type 0 vertex. Thus, we can denote the number of lazy walks with any such fixed signature $S$ by $F_p$.

We trivially observe that $F_0 = c$. Besides, we have $F_1 = \sum_{i = 1}^z b_i^2$ since there are $b_i^2$ ways to choose the starting and the ending type 1 vertex when the $i$-th type 0 vertex is selected to represent the zero block. Also, note that for any $1 \le i, j \le z$, the $i$-th and the $j$-th type 0 vertex have $\min\{ b_i, b_j \}$ common neighbors. Since the type 0 vertices representing distinct zero blocks can be selected independently, we conclude that
\begin{align}
    \label{fp_formula} F_p &= \sum_{i_1, i_2, \ldots, i_p = 1}^z b_{i_1} \min\{ b_{i_1}, b_{i_2} \} \min\{ b_{i_2}, b_{i_3} \} \cdots \min\{ b_{i_{p-1}}, b_{i_{p}} \} \, b_{i_p}\\
    \label{fp_formula_2} &= \sum_{i_1, i_2, \ldots, i_p = 1}^z b_{i_1} b_{\max\{i_1, i_2\}} b_{\max\{i_2, i_3\}} \cdots b_{\max\{i_{p-1}, i_p\}} \, b_{i_p}
\end{align}
holds for any $p \ge 2$. For more results on the expressions $F_p$, see Section~\ref{fp_stuff}.

We will now rely on the $(F_p)_{p \in \mathbb{N}_0}$ sequence to obtain a recurrence relation for computing the $\left(\LW_k\right)_{k \in \mathbb{N}_0}$ sequence. We proceed with the following auxiliary definitions.

\begin{Def}
    For any $k \in \mathbb{N}$, let $\Psi_k$ denote the set of all the lazy walks of length $k - 1$ in $G$ whose signature starts and ends with one.
\end{Def}
\begin{Def}
    For each $W = u_0 u_1 u_2 \cdots u_{k - 1} \in \Psi_k$, let $\sigma_1(W)$ be the longest lazy subwalk of $W$ of the form $u_0 u_1 u_2 \cdots u_{i - 1}$, where $1 \le i \le k$, which does not contain two consecutive type~1 vertices.
\end{Def}
\begin{Def}
    For each $W = u_0 u_1 u_2 \cdots u_{k - 1} \in \Psi_k$, where $\sigma_1(W) = u_0 u_1 u_2 \cdots u_{\ell - 1}$ for some $1 \le \ell \le k$, let $\sigma_2(W)$ be the lazy subwalk $u_{\ell} u_{\ell + 1} \cdots u_{k - 1}$ of $W$.
\end{Def}
\begin{Def}
    For any $k \in \mathbb{N}$ and $1 \le \ell \le k$, let $\Psi_{k, \ell}$ denote the subset of $\Psi_k$ comprising the lazy walks $W$ where $\sigma_1(W)$ is of length $\ell - 1$.
\end{Def}

For any $W = u_0 u_1 u_2 \cdots u_{k - 1} \in \Psi_{k, \ell}$, observe that if $\ell = k$, then $\sigma_2(W)$ is an empty lazy walk. On the other hand, if $\ell < k$, we have that $u_\ell$ is a type~1 vertex. Moreover, due to
\begin{equation}\label{lw_decomp}
    \LW_k = |\Psi_k| = \sum_{\ell = 1}^k |\Psi_{k, \ell}|,
\end{equation}
it is sufficient to count the walks from these $k$ subsets separately.

Since any two distinct type~1 vertices are adjacent, we obtain that for any fixed $1 \le \ell \le k$, the subwalk $\sigma_1(W)$ of a lazy walk $W = u_0 u_1 u_2 \cdots u_{k - 1}\in \Psi_{k, \ell}$ can be selected independently from the subwalk $\sigma_2(W)$. Moreover, $\sigma_2(W)$ can obviously be chosen in $\LW_{k - \ell}$ different ways. Note that $\sigma_2(W)$ is empty when $\ell = k$, which is consistent with $\LW_0 = 1$. On the other hand, the signature of $\sigma_1(W)$ is of the form \eqref{f_form} with a certain number of zero blocks. If $\ell = 1$, there are $F_0 = c$ choices for $\sigma_1(W)$. For any $\ell \ge 2$, it can be routinely checked that there are $\binom{\ell - p - 2}{p - 1}$ signatures of the form~\eqref{f_form} with $\ell$ vertices and $p \in \mathbb{N}$ zero blocks, with each such signature corresponding to $F_p$ lazy walks. Thus, we get $|\Psi_{k, 1}| = \LW_{k - 1} \, F_0$, alongside
\[
    |\Psi_{k, \ell}| = \LW_{k - \ell} \sum_{p \in \mathbb{N}} \binom{\ell - p - 2}{p - 1} F_p \qquad (2 \le \ell \le k) .
\]
Observe that for any $k \ge 2$, we have $|\Psi_{k, 2}| = 0$ because $\sigma_1(W)$ cannot consist of two consecutive type 1 vertices. With all of this in mind, \eqref{lw_decomp} gives us
\[
    \LW_k = \LW_{k-1} \, F_0 + \sum_{\ell = 3}^k \LW_{k - \ell} \sum_{p \in \mathbb{N}} \binom{\ell - p - 2}{p - 1} F_p = \LW_{k-1} \, F_0 + \sum_{r = 0}^{k - 3} \LW_r \sum_{q \in \mathbb{N}_0} \binom{k - 3 - r - q}{q} F_{q + 1} .
\]
From here, we obtain the following.

\begin{Prop}\label{recurrence_prop}
The $\left(\LW_k\right)_{k \in \mathbb{N}_0}$ sequence can be computed by the recurrence relation
\begin{align}
    \label{lw_rec_1} \LW_0 &= 1,\\
    \label{lw_rec_2} \LW_k &= c \, \LW_{k - 1} + \sum_{r = 0}^{k - 3} \LW_r \sum_{q \in \mathbb{N}_0} \binom{k - 3 - r - q}{q} F_{q + 1} \qquad (k \in \mathbb{N}).
\end{align}
\end{Prop}

\section{Lower bound via lazy walks}\label{sc_lower_bound}

In this section, our goal will be to get a lower bound for $\rho$ through the next theorem.

\begin{Th}\label{cubic_th_1}
    The greatest real root of the polynomial
    \begin{equation}\label{cubic_poly_1}
        x^3 - (c+1)x^2 + cx - F_1
    \end{equation}
    is a lower bound for $1 + \rho$.
\end{Th}

To begin, let $\left(\LW'_k\right)_{k \in \mathbb{N}_0}$ be the auxiliary sequence defined by
\begin{alignat}{2}
    \label{lw1_rec_1}\LW'_k &= c^k && \qquad (0 \le k \le 2),\\
    \label{lw1_rec_2}\LW'_k &= c \, \LW'_{k - 1} + \sum_{r = 0}^{k - 3} \LW'_r \, F_1 && \qquad (k \ge 3).
\end{alignat}
By comparing \eqref{lw1_rec_1} and \eqref{lw1_rec_2} with \eqref{lw_rec_1} and \eqref{lw_rec_2} from Proposition \ref{recurrence_prop}, it can be routinely verified by induction that
\begin{equation}\label{lw_boundedness_1}
    \LW'_k \le \LW_k \qquad (k \in \mathbb{N}_0) .
\end{equation}
Also, we have the following lemma.

\begin{Lemma}\label{cubic_rec_1}
The $\left(\LW'_k\right)_{k \in \mathbb{N}_0}$ sequence satisfies the linear recurrence relation
\begin{equation}\label{lw1_rec_3}
    \LW'_k - (c + 1) \, \LW'_{k - 1} + c \, \LW'_{k - 2} - F_1 \, \LW'_{k - 3} = 0 \qquad (k \ge 3).
\end{equation}
\end{Lemma}
\begin{proof}
Due to \eqref{lw1_rec_2}, for every $k \ge 3$ we have
\begin{align*}
    \LW'_k - \LW'_{k - 1} &= \left( c \, \LW'_{k - 1} + F_1 \sum_{r = 0}^{k - 3} \LW'_r \right) - \left( c \, \LW'_{k - 2} + F_1 \sum_{r = 0}^{k - 4} \LW'_r \right)\\
    &= c \, \LW'_{k - 1} - c \, \LW'_{k - 2} + F_1 \, \LW'_{k - 3} .
\end{align*}
Thus, \eqref{lw1_rec_3} follows directly from here.
\end{proof}

We resume by obtaining the next result concerning the characteristic polynomial of the recurrence relation \eqref{lw1_rec_3}.

\begin{Lemma}\label{cubic_lemma_1}
    Let $\xi \in \mathbb{R}$ be the greatest real root of \eqref{cubic_poly_1}. Then $\xi$ is a simple root of \eqref{cubic_poly_1} and has a greater modulus than all the other roots.
\end{Lemma}
\begin{proof}
    Let
    \[
        P(x) = x^3 - (c+1)x^2 + cx - F_1 = x(x-1)(x-c) - F_1.
    \]
    Since $P(c) = -F_1 < 0$ and the function $x \mapsto P(x)$ is increasing on $[c, +\infty)$, it follows that $\xi > c$. The formal derivative
    \[
        P'(x) = 3x^2 - 2(c+1)x + c = 3\left( x - \frac{c+1}{3}\right)^2 + \frac{-c^2+c-1}{3}
    \]
    is such that $P'(x) \ge P'(c) = c(c-1) > 0$ for every $x \ge c$. Thus, $P'(\xi) > 0$, which implies that $\xi$ is a simple root of $P(x)$.

    Let the other two roots of $P(x)$ be $\xi_2$ and $\xi_3$. Since $P(x) < 0$ for every $x \in (-\infty, 0]$, we have that if $\xi_2, \xi_3 \in \mathbb{R}$, then these two roots are from $(0, c)$, hence they have a smaller modulus than $\xi$. On the other hand, if $\xi_2, \xi_3 \in \mathbb{C} \setminus \mathbb{R}$, then these two numbers are a complex conjugate pair. Now, suppose that $|\xi_2| = |\xi_3| \ge \xi$. In this case, Vieta's formulas yield
    \[
        F_1 = \xi \xi_2 \xi_3 = \xi \, |\xi_2|^2 \ge \xi^3,
    \]
    which further gives
    \[
        0 = P(\xi) = \xi^3 - (c+1)\xi^2 + c \, \xi - F_1 = (\xi^3 - F_1) + \xi\left( c - (c+1) \xi \right) \le \xi\left( c - (c+1) \xi \right) .
    \]
    Therefore, $c - (c + 1) \xi \ge 0$, which is impossible due to $\xi > c \ge 2$. This contradiction guarantees that $\xi$ must have a greater modulus than all the other roots of $P(x)$.
\end{proof}

We are now in a position to finalize the proof of Theorem \ref{cubic_th_1}.

\begin{proof}[Proof of Theorem~\ref{cubic_th_1}]
    Let $\xi$ be the greatest real root of \eqref{cubic_poly_1} and let the remaining two roots be $\xi_2$ and $\xi_3$. From Lemmas~\ref{cubic_rec_1} and \ref{cubic_lemma_1}, we obtain either
    \[
        \mathrm{LW}'_k = \alpha_1 \xi^k + \alpha_2 \xi_2^k + \alpha_3 \xi_3^k \qquad (k \in \mathbb{N}_0),
    \]
    if $\xi_2 \neq \xi_3$, or
    \[
        \mathrm{LW}'_k = \alpha_1 \xi^k + \xi_2^k (\alpha_2 + \alpha_3 k) \qquad (k \in \mathbb{N}_0),
    \]
    if $\xi_2 = \xi_3$, for some $\alpha_1, \alpha_2, \alpha_3 \in \mathbb{C}$. Either way, we claim that $\alpha_1 \neq 0$. Indeed, if we had $\alpha_1 = 0$, then $\left(\LW'_k\right)_{k \in \mathbb{N}_0}$ would be a solution to some linear recurrence relation of order two with the characteristic polynomial
    \[
        x^2 + \beta x + \gamma,
    \]
    where $\beta, \gamma \in \mathbb{R}$. From \eqref{lw1_rec_1}, we would get $c^2 + \beta \, c + \gamma = 0$, thus giving
    \[
        c^k + \beta \, c^{k - 1} + \gamma \, c^{k - 2} = 0 \qquad (k \ge 2) .
    \]
    Therefore, the sequence $\left(c^k\right)_{k \in \mathbb{N}_0}$ would also be a solution to this linear recurrence relation of order two at the same time as $\left(\LW'_k\right)_{k \in \mathbb{N}_0}$. This is impossible due to $\LW'_3 = c^3 + F_1 \neq c^3$, since there exists a unique solution with two predetermined starting terms. Thus, we obtain $\alpha_1 \neq 0$. Taking everything into consideration, the approach from Lemma \ref{lazy_lemma} can analogously be used to reach
    \[
        \lim_{k \to +\infty} \sqrt[k]{\LW'_k} = \xi.
    \]
    The result now follows directly from Lemma \ref{lazy_lemma} and \eqref{lw_boundedness_1}.
\end{proof}

We end the section with the next corollary.

\begin{Cor}
    We have
    \[
        \rho > c - 1 + \frac{F_1}{n^2} .
    \]
\end{Cor}
\begin{proof}
    Let
    \[
        P(x) = x^3 - (c+1)x^2 + cx - F_1 = x(x-1)(x-c) - F_1
    \]
    and $x_0 = c + \frac{F_1}{n^2}$. Due to Theorem \ref{cubic_th_1}, it is enough to prove that $x_0$ is smaller than the greatest real root of $P(x)$. Thus, it suffices to show that $P(x_0) < 0$. Since
    \[
        F_1 = \sum_{i = 1}^z b_i^2 < (n - c) c^2 < (n - c) n^2 ,
    \]
    we obtain $\frac{F_1}{n^2} < n - c$, i.e., $x_0 < n$. This leads us to
    \[
        \pushQED{\qed}
        P(x_0) = (x_0^2 - x_0)(x_0 - c) - F_1 < x_0^2 (x_0 - c) - F_1 < n^2 \cdot \frac{F_1}{n^2} - F_1 = 0.\qedhere
    \]
\end{proof}

\section{Upper bound via lazy walks}\label{sc_upper_bound}

In the present section, we will obtain the following upper bound for $\rho$.

\begin{Th}\label{cubic_th_2}
    The greatest real root of the polynomial
    \begin{equation}\label{cubic_poly_2}
        x^3 - (c+1)x^2 + \left( c - \sum_{i = 1}^z b_i \right) x + \left( c \sum_{i = 1}^z b_i - F_1 \right)
    \end{equation}
    is an upper bound for $1 + \rho$.
\end{Th}

We start by defining the auxiliary sequence $\left(\LW''_k\right)_{k \in \mathbb{N}_0}$ via
\begin{alignat}{2}
    \label{lw2_rec_1} \LW''_0 &= 1,\\
    \label{lw2_rec_2} \LW''_k &= c \, \LW''_{k - 1} + \sum_{r = 0}^{k - 3} \LW''_r \sum_{q \in \mathbb{N}_0} \binom{k - 3 - r - q}{q} F_1 \left( \sum_{i = 1}^z b_i \right)^{q} \qquad (k \in \mathbb{N}).
\end{alignat}
This sequence is actually an upper bound for $\left(\LW_k\right)_{k \in \mathbb{N}_0}$, as shown in the following lemma.

\begin{Lemma}\label{cubic_ineq_2}
    For every $k \in \mathbb{N}_0$, we have $\LW''_k \ge \LW_k$.
\end{Lemma}
\begin{proof}
    Since $\min\{b_i, b_j\} \le \sqrt{b_i b_j}$ for any $1 \le i, j \le z$, \eqref{fp_formula} further gives us
    \begin{align*}
        F_p &\le \sum_{i_1, i_2, \ldots, i_p = 1}^z b_{i_1} \sqrt{b_{i_1} b_{i_2}} \sqrt{b_{i_2} b_{i_3}} \cdots \sqrt{b_{i_{p-1}} b_{i_p}} \,  b_{i_p} = \sum_{i_1, i_2, \ldots, i_p = 1}^z b_{i_1}^\frac{3}{2} b_{i_2} b_{i_3} \cdots b_{i_{p-1}} b_{i_p}^\frac{3}{2}\\
        &= \left( \sum_{i = 1}^z b_i^\frac{3}{2} \right)^2 \left( \sum_{i = 1}^z b_i \right)^{p - 2} = \left( \sum_{i = 1}^z b_i \sqrt{b_i} \right)^2 \left( \sum_{i = 1}^z b_i \right)^{p - 2}
    \end{align*}
    for every $p \ge 2$. Also, the Cauchy--Schwarz inequality yields
    \[
        \left( \sum_{i = 1}^z b_i \sqrt{b_i} \right)^2 \le \left( \sum_{i = 1}^z b_i^2 \right) \left( \sum_{i = 1}^z b_i \right) = F_1 \left( \sum_{i = 1}^z b_i \right) .
    \]
    Thus, we obtain
    \begin{equation}\label{fp_cauchy}
        F_p \le F_1 \left( \sum_{i = 1}^z b_i \right)^{p - 1} \qquad (p \ge 2).
    \end{equation}
    The result can now be routinely verified by induction using \eqref{fp_cauchy} and by comparing \eqref{lw2_rec_1} and \eqref{lw2_rec_2} with \eqref{lw_rec_1} and \eqref{lw_rec_2} from Proposition~\ref{recurrence_prop}.
\end{proof}

We now need the next result.

\begin{Lemma}\label{cubic_rec_2}
The sequence $\left(\LW''_k\right)_{k \in \mathbb{N}_0}$ satisfies the linear recurrence relation
\begin{equation}\label{lw2_rec_3}
    \LW''_k - (c + 1) \, \LW''_{k - 1} + \left( c - \sum_{i = 1}^z b_i \right) \LW''_{k - 2} + \left( c \sum_{i = 1}^z b_i - F_1 \right) \LW''_{k - 3} = 0 \qquad (k \ge 3).
\end{equation}
\end{Lemma}
\begin{proof}
Due to \eqref{lw2_rec_2}, for every $k \ge 3$ we have
\begin{align*}
    \LW''_k - \LW''_{k - 1} &= c \, \LW''_{k - 1} + F_1 \sum_{r = 0}^{k - 3} \LW''_r \sum_{q \in \mathbb{N}_0} \binom{k - 3 - r - q}{q} \left( \sum_{i = 1}^z b_i \right)^{q}\\
    &- c \, \LW''_{k - 2} - F_1 \sum_{r = 0}^{k - 4} \LW''_r \sum_{q \in \mathbb{N}_0} \binom{k - 4 - r - q}{q} \left( \sum_{i = 1}^z b_i \right)^{q} .
\end{align*}
This leads us to
\begin{align*}
    \LW''_k - (c + 1) \, \LW''_{k - 1} &+ c \, \LW''_{k - 2} = F_1 \, \LW''_{k - 3} \sum_{q \in \mathbb{N}_0} \binom{-q}{q} \left( \sum_{i = 1}^z b_i \right)^{q}\\
    &+ F_1 \sum_{r = 0}^{k - 4} \LW''_r \sum_{q \in \mathbb{N}_0} \left( \binom{k - 3 - r - q}{q} - \binom{k - 4 - r - q}{q} \right) \left( \sum_{i = 1}^z b_i \right)^{q},
\end{align*}
i.e.,
\begin{align}
    \nonumber \LW''_k - (c + 1) \, \LW''_{k - 1} &+ c \, \LW''_{k - 2} - F_1 \, \LW''_{k - 3} = F_1 \sum_{r = 0}^{k - 4} \LW''_r \sum_{q \in \mathbb{N}} \binom{k - 4 - r - q}{q - 1} \left( \sum_{i = 1}^z b_i \right)^{q}\\
    \nonumber &= F_1 \sum_{r = 0}^{k - 4} \LW''_r \sum_{q' \in \mathbb{N}_0} \binom{k - 5 - r - q'}{q'} \left( \sum_{i = 1}^z b_i \right)^{q' + 1}\\
    \label{lw2_weird_1} &= \left( \sum_{i = 1}^z b_i \right) F_1 \sum_{r = 0}^{k - 5} \LW''_r \sum_{q' \in \mathbb{N}_0} \binom{k - 5 - r - q'}{q'} \left( \sum_{i = 1}^z b_i \right)^{q'} .
\end{align}
At the same time, \eqref{lw2_rec_2} directly yields
\begin{equation}\label{lw2_weird_2}
    \LW''_{k - 2} - c \, \LW''_{k - 3} = F_1 \sum_{r = 0}^{k - 5} \LW''_r \sum_{q \in \mathbb{N}_0} \binom{k - 5 - r - q}{q} \left( \sum_{i = 1}^z b_i \right)^{q}.
\end{equation}
By combining \eqref{lw2_weird_1} and \eqref{lw2_weird_2}, we conclude that
\[
    \LW''_k - (c + 1) \, \LW''_{k - 1} + c \, \LW''_{k - 2} - F_1 \, \LW''_{k - 3} = \left( \sum_{i = 1}^z b_i \right) \left( \LW''_{k - 2} - c \, \LW''_{k - 3} \right)
\]
holds for any $k \ge 3$. Therefore, \eqref{lw2_rec_3} follows immediately from here.
\end{proof}

We proceed by deriving the following result regarding the characteristic polynomial of the recurrence relation~\eqref{lw2_rec_3}.

\begin{Lemma}\label{cubic_lemma_2}
    Let $\xi \in \mathbb{R}$ be the greatest real root of \eqref{cubic_poly_2}. Then $\xi$ is a simple root of \eqref{cubic_poly_2} and has a greater modulus than all the other roots.
\end{Lemma}
\begin{proof}
    Let
    \[
        P(x) = x^3 - (c+1)x^2 + \left( c - \sum_{i = 1}^z b_i \right) x + \left( c \sum_{i = 1}^z b_i - F_1 \right).
    \]
    Observe that $P(c) = -F_1 < 0$, as well as
    \[
        P(0) = c \sum_{i = 1}^z b_i - F_1 = c \sum_{i = 1}^z b_i - \sum_{i = 1}^z b_i^2 = \sum_{i = 1}^z b_i (c - b_i) > 0.
    \]
    It follows from here that $P(x)$ has three simple real roots $\xi$, $\xi_2$, $\xi_3$, such that $\xi \in (c, +\infty)$, $\xi_2 \in (0, c)$ and $\xi_3 \in (-\infty, 0)$. Obviously, we have $|\xi_2| < |\xi|$. Now, suppose that $|\xi_3| \ge |\xi|$. In this case, we get $-\xi_3 \ge \xi$, i.e., $\xi + \xi_3 \le 0$. Since $\xi_2 < c$, this means that $\xi + \xi_2 + \xi_3 < c$, which is impossible because Vieta's formulas directly yield $\xi + \xi_2 + \xi_3 = c + 1$. Therefore, the real root $\xi$ must have a greater modulus than all the other roots.
\end{proof}

Theorem \ref{cubic_th_2} can now be proved analogously to Theorem \ref{cubic_th_1}, so we omit the details of its proof.

\section{Two additional lower bounds}\label{sc_extra_bounds}
In the last section, we provide one more lower bound for $\rho$, as well as an alternative proof of Theorem \ref{cubic_th_2}. Here, we will assume that $V(G) = \{1, 2, \ldots, n \}$ and consider the vertices to be arranged in descending order with respect to their degree so that $n - 1 = d_1 \ge d_2 \ge \cdots \ge d_n \ge 1$ and all the type 1 vertices appear before the type 0 vertices. Observe that in this case, the obtained adjacency matrix $A \in \mathbb{R}^{n \times n}$ is of stepwise form. Hence, we have $d_c = c - 1$, together with $d_{c + i} = b_i$ for any $1 \le i \le z$. Moreover, the corresponding Perron--Frobenius vector $w \in \mathbb{R}^n$ is such that $w_1 \ge w_2 \ge \cdots \ge w_n > 0$ (see, e.g., \cite{brusol1986}). With all of this in mind, we obtain the following two new lower bounds.

\begin{Th}\label{last_th}
We have
\begin{equation}\label{quad_lower_formula}
    \rho \ge \frac{c - 2 + \sqrt{c^2 + \frac{4}{c - 1} F_1}}{2}.
\end{equation}
\end{Th}
\begin{proof}
    For the Perron--Frobenius vector $w \in \mathbb{R}^n$, we have
    \begin{alignat}{2}
        \label{perron_1} (\rho + 1) w_i &= \sum_{j = 1}^{d_i + 1} w_j \qquad && (1 \le i \le c - 1),\\
        \label{perron_2} \rho w_i &= \sum_{j = 1}^{d_i} w_j \qquad && (c \le i \le n) .
    \end{alignat}
    By summing up \eqref{perron_1} for $i = 1, 2, \ldots, c - 1$, we further get
    \begin{equation}\label{cool_1}
        (\rho + 1) \sum_{i = 1}^{c - 1} w_i = (c - 1) \sum_{i = 1}^{c - 1} w_i + \sum_{i = c}^n d_i w_i .
    \end{equation}
    Also, since $w_1 \ge w_2 \ge \cdots \ge w_n > 0$ and $c - 1 = d_c \ge d_{c + 1} \ge \cdots \ge d_n$, it follows that
    \[
        \frac{1}{d_i} \sum_{j = 1}^{d_i} w_j \ge \frac{1}{c - 1} \sum_{j = 1}^{c - 1} w_j
    \]
    holds for any $c \le i \le n$. Therefore, \eqref{perron_2} implies
    \begin{equation}\label{cool_2}
        w_i = \frac{1}{\rho} \sum_{j = 1}^{d_i} w_j \ge \frac{d_i}{\rho(c - 1)} \sum_{j = 1}^{c - 1} w_j
    \end{equation}
    for every $c \le i \le n$.

    By combining \eqref{cool_1} and \eqref{cool_2}, we obtain
    \[
        (\rho + 1) \sum_{i = 1}^{c - 1} w_i \ge (c - 1) \sum_{i = 1}^{c - 1} w_i + \sum_{i = c}^n \frac{d_i^2}{\rho(c-1)} \sum_{j = 1}^{c - 1} w_j ,
    \]
    which leads us to
    \[
        \rho + 1 \ge c - 1 + \sum_{i = c}^n \frac{d_i^2}{\rho(c-1)} .
    \]
    Note that $F_1 = \sum_{i = 1}^z b_i^2 = \sum_{i = c + 1}^n d_i^2$ and $d_c^2 = (c-1)^2$. Since $\rho > 0$, we further have
    \[
        \rho(\rho + 1) \ge \rho (c - 1) + \frac{F_1 + (c - 1)^2}{c - 1},
    \]
    i.e.,
    \begin{equation}\label{aux_bound}
        \rho^2 - (c - 2)\rho - \left( \frac{F_1}{c - 1} + c - 1 \right) \ge 0 .
    \end{equation}
    The polynomial $x^2 - (c - 2)x - \frac{F_1 + (c - 1)^2}{c - 1}$ obviously has a positive and a negative simple root, which means that \eqref{aux_bound} is equivalent to \eqref{quad_lower_formula}.
\end{proof}

\begin{Th}
We have
\begin{equation}\label{ugly_bound}
    \rho\left( (\rho - c + 2)(\rho^2 + \rho - (z+1)) - S \right)(c - 2) \le \sum_{i = c}^n \left( d_i (\rho^2 - (z+1)) - \rho(\rho - c + 2) + S \right) (d_i - 1) ,
\end{equation}
where $S = \sum_{i = c}^n d_i = c - 1 + \sum_{i = 1}^z b_i$.
\end{Th}
\begin{proof}
    Let $w \in \mathbb{R}^n$ be the Perron--Frobenius vector, and let $A_1 = \sum_{i = 1}^{c - 1} w_i$ and $A_2 = \sum_{i = 2}^{c - 1} w_i$. In the same way as in Theorem \ref{last_th}, we may obtain \eqref{cool_1}, which then leads us to
    \begin{equation}\label{aux_1}
        \rho (\rho - c + 2) A_1 = \sum_{i = c}^n \rho d_i w_i.
    \end{equation}
    By plugging in~\eqref{perron_2} into \eqref{aux_1}, we obtain
    \begin{equation}\label{aux_2}
        \rho(\rho - c + 2) A_1 = \sum_{i = c}^n d_i \sum_{j = 1}^{d_i} w_j = \sum_{i = c}^n d_i \left( w_1 + \sum_{j = 2}^{d_i} w_j \right) = S w_1 + \sum_{i = c}^n d_i \sum_{j = 2}^{d_i} w_j .
    \end{equation}
    Since $A_1 = A_2 + w_1$, \eqref{aux_2} transforms into
    \begin{equation}\label{aux_3}
        \rho (\rho - c + 2)A_2 + \left( \rho(\rho - c + 2) - S \right) w_1 = \sum_{i = c}^n d_i \sum_{j = 2}^{d_i} w_j .
    \end{equation}

    The vertex $1$ is adjacent to all the other vertices, hence we have
    \[
        \rho w_1 = A_2 + \sum_{i = c}^n w_i ,
    \]
    which together with \eqref{perron_2} further gives
    \[
        \rho^2 w_1 = \rho A_2 + \sum_{i = c}^n \rho w_i = \rho A_2 + \sum_{i = c}^n \sum_{j = 1}^{d_i} w_j = \rho A_2 + \sum_{i = c}^n \left( w_1 + \sum_{j = 2}^{d_i} w_j \right) = \rho A_2 + (z + 1) w_1 + \sum_{i = c}^n \sum_{j = 2}^{d_i} w_j .
    \]
    Thus, we obtain
    \begin{equation}\label{aux_4}
        (\rho^2 - (z + 1)) w_1 = \rho A_2 + \sum_{i = c}^n \sum_{j = 2}^{d_i} w_j .
    \end{equation}
    Note that the right-hand side of \eqref{aux_4} is positive, which assures us that $\rho^2 - (z + 1)$ is also positive.

    By multiplying \eqref{aux_3} with $\rho^2 - (z + 1)$ and plugging in \eqref{aux_4}, we reach
    \[
        \rho(\rho - c + 2)(\rho^2 - (z + 1)) A_2 + \left(\rho(\rho - c + 2) - S \right) \left( \rho A_2 + \sum_{i = c}^n \sum_{j = 2}^{d_i} w_j \right) = \sum_{i = c}^n d_i (\rho^2 - (z + 1)) \sum_{j = 2}^{d_i} w_j,
    \]
    i.e.,
    \begin{equation}\label{aux_5}
       \rho\left( (\rho - c + 2)(\rho^2 + \rho - (z+1)) - S \right) A_2 = \sum_{i = c}^n \left( d_i (\rho^2 - (z+1)) - \rho(\rho - c + 2) + S \right) \sum_{j = 2}^{d_i} w_j .
    \end{equation}
    Observe that for any $c \le i \le n$, we have
    \begin{align*}
        d_i (\rho^2 - (z+1)) - \rho(\rho - c + 2) + S &\ge (\rho^2 - (z + 1)) - \rho(\rho - c + 2) + S\\
        &= \rho(c - 2) + (S - (z + 1)) \ge \rho (c - 2) > 0.
    \end{align*}
    Therefore, we may plug in the inequality
    \[
        \sum_{j = 2}^{d_i} w_j \ge \frac{d_i - 1}{c - 2} \sum_{j = 2}^{c - 1} w_j = \frac{d_i - 1}{c - 2} A_2 \qquad (c \le i \le n)
    \]
    into \eqref{aux_5} in order to obtain
    \[
       \rho\left( (\rho - c + 2)(\rho^2 + \rho - (z+1)) - S \right) A_2 = \sum_{i = c}^n \left( d_i (\rho^2 - (z+1)) - \rho(\rho - c + 2) + S \right) \frac{d_i - 1}{c - 2} A_2 .
    \]
    From here, we get \eqref{ugly_bound}.
\end{proof}

We end the section by giving an alternative proof of Theorem \ref{cubic_th_2}.

\begin{proof}[Alternative proof of Theorem \ref{cubic_th_2}]
    Let
    \[
        P(x) = x^3 - (c+1)x^2 + \left( c - \sum_{i = 1}^z b_i \right) x + \left( c \sum_{i = 1}^z b_i - F_1 \right) = (x - c)\left( x(x-1) - \sum_{i = 1}^z b_i \right) - F_1 .
    \]
    Clearly, it is sufficient to show that $P(\rho + 1) \le 0$, i.e.,
    \begin{equation}\label{alt_bound}
        (\rho - c + 1)\left( \rho(\rho + 1) - \sum_{i = 1}^z b_i \right) \le F_1 .
    \end{equation}
    For the Perron--Frobenius vector $w \in \mathbb{R}$, we can reach \eqref{cool_1} in the same way as in Theorem~\ref{last_th}. Moreover, for any $c + 1 \le i \le n$, \eqref{perron_2} implies
    \begin{equation}\label{cool_3}
        \rho w_i = \sum_{j = 1}^{d_i} w_j \le d_i w_1 .
    \end{equation}
    By combining \eqref{cool_1} with \eqref{cool_3}, we obtain
    \[
        (\rho + 1) \sum_{i = 1}^c w_i \le c \sum_{i = 1}^c w_i + \sum_{i = c + 1}^n \frac{d_i^2}{\rho} w_1 ,
    \]
    i.e.,
    \begin{equation}\label{cool_4}
        (\rho - c + 1) \sum_{i = 1}^c w_i \le \frac{F_1}{\rho} w_1 .
    \end{equation}

    Since $d_1 = n - 1$, \eqref{perron_1} gives
    \begin{equation}\label{cool_5}
        (\rho + 1) w_1 = \sum_{j = 1}^n w_j = \sum_{j = 1}^c w_j + \sum_{i = c + 1}^n w_i .
    \end{equation}
    By combining \eqref{cool_5} with \eqref{cool_3}, we get
    \[
        (\rho + 1) w_1 \le \sum_{j = 1}^c w_j + \sum_{i = c + 1}^{n} \frac{d_i}{\rho} w_1 .
    \]
    Since $d_{c + i} = b_i$ for any $1 \le i \le z$, it follows that
    \begin{equation}\label{cool_6}
        \left( \rho + 1 - \frac{1}{\rho} \sum_{i = 1}^{z} b_i \right) w_1 \le \sum_{j = 1}^c w_j .
    \end{equation}
    
    We now finalize the proof by splitting the problem into two cases as follows. If $\rho + 1 - \frac{1}{\rho} \sum_{i = 1}^{z} b_i \le 0$, we have $\rho(\rho + 1) - \sum_{i = 1}^z b_i \le 0$, hence to prove \eqref{alt_bound}, it is enough to show that $\rho \ge c - 1$. This directly follows from the fact that $G$ contains $K_c$ as a subgraph (see, e.g., \cite[Chapter 8]{godroy2001}). On the other hand, if $\rho + 1 - \frac{1}{\rho} \sum_{i = 1}^{z} b_i > 0$, then \eqref{cool_6} becomes
    \begin{equation}\label{cool_7}
        w_1 \le \frac{\rho}{\rho(\rho + 1) - \sum_{i = 1}^z b_i} \sum_{j = 1}^c w_j .
    \end{equation}
    By plugging in \eqref{cool_7} into \eqref{cool_4}, we obtain
    \[
        \rho - c + 1 \le \frac{F_1}{\rho} \cdot \frac{\rho}{\rho(\rho + 1) - \sum_{i = 1}^z b_i} ,
    \]
    which immediately yields \eqref{alt_bound}.
\end{proof}

\appendix
\section{\texorpdfstring{$F_p$}{F(p)} expressions and FOP sequence}\label{fp_stuff}

In the present section, we provide a few results concerning the computation of the $F_p$ expressions. To begin, let $B \in \mathbb{R}^{z \times z}$ be the matrix defined by
\[
    B_{ij} = b_{\max \{i, j\}} \qquad (1 \le i, j \le z) ,
\]
so that
\[
    B = \begin{bmatrix}
        b_1 & b_2 & b_3 & b_4 & \cdots & b_z\\
        b_2 & b_2 & b_3 & b_4 & \cdots & b_z\\
        b_3 & b_3 & b_3 & b_4 & \cdots & b_z\\
        b_4 & b_4 & b_4 & b_4 & \cdots & b_z\\
        \vdots & \vdots & \vdots & \vdots & \ddots & \vdots\\
        b_z & b_z & b_z & b_z & \cdots & b_z
    \end{bmatrix} .
\]
With this in mind, \eqref{fp_formula_2} immediately implies the next proposition.

\begin{Prop}
    For any $p \in \mathbb{N}$, we have
    \[
        F_p = w^\intercal \, B^{p - 1} \, w,
    \]
    where $w = \begin{bmatrix} b_1 & b_2 & b_3 & \cdots & b_z \end{bmatrix}^\intercal$.
\end{Prop}

From the spectral decomposition of $B$, we also obtain the following corollary.

\begin{Cor}
    Let $\lambda_1, \lambda_2, \ldots, \lambda_z$ be the eigenvalues of $B$ with the corresponding eigenvectors $x_1, x_2, \ldots, x_z$, and let $w = \begin{bmatrix} b_1 & b_2 & b_3 & \cdots & b_z \end{bmatrix}^\intercal$. Then we have
    \[
        F_p = \sum_{i = 1}^z (w^\intercal x_i^\intercal x_i w) \, \lambda_i^{p - 1} \qquad (p \in \mathbb{N}).
    \]
\end{Cor}

Observe that $B$ is positive semidefinite. Indeed, we have
\begin{align}\label{sumb}
\begin{split}
    B = \begin{bmatrix}
        b_z & b_z & \cdots & b_z & b_z\\
        b_z & b_z & \cdots & b_z & b_z\\
        \vdots & \vdots & \ddots & \vdots & \vdots\\
        b_z & b_z & \cdots & b_z & b_z\\
        b_z & b_z & \cdots & b_z & b_z
    \end{bmatrix} &+ \begin{bmatrix}
        b_{z - 1} - b_z & b_{z - 1} - b_z & \cdots & b_{z - 1} - b_z & 0\\
        b_{z - 1} - b_z & b_{z - 1} - b_z & \cdots & b_{z - 1} - b_z & 0\\
        \vdots & \vdots & \ddots & \vdots & \vdots\\
        b_{z - 1} - b_z & b_{z - 1} - b_z & \cdots & b_{z - 1} - b_z & 0\\
        0 & 0 & \cdots & 0 & 0
    \end{bmatrix} + \cdots\\
    &+ \begin{bmatrix}
        b_2 - b_3 & b_2 - b_3 & 0 & \cdots & 0\\
        b_2 - b_3 & b_2 - b_3 & 0 & \cdots & 0\\
        0 & 0 & 0 & \cdots & 0\\
        \vdots & \vdots & \vdots & \ddots & \vdots\\
        0 & 0 & 0 & \cdots & 0
    \end{bmatrix} + \begin{bmatrix}
        b_1 - b_2 & 0 & 0 & \cdots & 0\\
        0 & 0 & 0 & \cdots & 0\\
        0 & 0 & 0 & \cdots & 0\\
        \vdots & \vdots & \vdots & \ddots & \vdots\\
        0 & 0 & 0 & \cdots & 0
    \end{bmatrix} .
\end{split}
\end{align}
Since $b_1 \ge b_2 \ge \cdots \ge b_z$, it follows that each summand matrix from \eqref{sumb} is positive semidefinite, hence $B$ itself is also positive semidefinite.

We now elaborate another way how the $(F_p)_{p \in \mathbb{N}_0}$ elements can be expressed. Analogously to the BZP sequence, we introduce the \emph{forward one position} (FOP) sequence as the nondecreasing tuple $(f_1, f_2, \ldots, f_c) \in \mathbb{N}^c$ where $f_i$ denotes the number of type~0 vertices that appear before the $i$-th type~1 vertex in the generating sequence. Note that $f_1 = 0$ and $f_c = z$. Clearly, the generating sequence can always be reconstructed from the given FOP sequence and the number of vertices. Moreover, we have $\sum_{i = 1}^c f_i = \sum_{i = 1}^z b_i = m - \binom{c}{2}$.

The number of lazy walks whose signature is of the form \eqref{f_form} can now be counted alternatively as follows. For each choice of $p + 1$ type~1 vertices, the type~0 vertices corresponding to the $p$ zero blocks can be selected independently. Furthermore, for any $1 \le i, j \le c$, the $i$-th and the $j$-th type~1 vertex have $\min\{ f_i, f_j \}$ common neighbors among the type~0 vertices. Thus, we obtain
\begin{align}
    \nonumber F_p &= \sum_{i_1, i_2, \ldots, i_{p + 1} = 1}^c \min\{ f_{i_1}, f_{i_2} \} \min\{ f_{i_2}, f_{i_3} \} \cdots \min\{ f_{i_p}, f_{i_{p + 1}} \}\\
    \label{fop_thing} &= \sum_{i_1, i_2, \ldots, i_{p + 1} = 1}^c f_{\min\{i_1, i_2\}} f_{\min\{i_2, i_3\}} \cdots f_{\min\{i_p, i_{p+1}\}}
\end{align}
for every $p \in \mathbb{N}$.

Let $\Phi \in \mathbb{R}^{c \times c}$ be the matrix defined by
\[
    \Phi_{ij} = f_{\min \{i, j\}} \qquad (1 \le i, j \le c) ,
\]
so that
\[
    \Phi = \begin{bmatrix}
        f_1 & f_1 & f_1 & f_1 & \cdots & f_1\\
        f_1 & f_2 & f_2 & f_2 & \cdots & f_2\\
        f_1 & f_2 & f_3 & f_3 & \cdots & f_3\\
        f_1 & f_2 & f_3 & f_4 & \cdots & f_4\\
        \vdots & \vdots & \vdots & \vdots & \ddots & \vdots\\
        f_1 & f_2 & f_3 & f_4 & \cdots & f_c
    \end{bmatrix} .
\]
Note that $\Phi$ is also positive semidefinite. From~\eqref{fop_thing} we may analogously get the following proposition and corollary.

\begin{Prop}
    For any $p \in \mathbb{N}_0$, we have
    \[
        F_p = \bm{1}^\intercal \, \Phi^p \, \bm{1},
    \]
    where $\bm{1} = \begin{bmatrix} 1 & 1 & 1 & \cdots & 1 \end{bmatrix}^\intercal \in \mathbb{R}^{c \times 1}$.
\end{Prop}
\begin{Cor}
    Let $\lambda_1, \lambda_2, \ldots, \lambda_c$ be the eigenvalues of $\Phi$ with the corresponding eigenvectors $x_1, x_2, \ldots, x_c$, and let $\bm{1} = \begin{bmatrix} 1 & 1 & 1 & \cdots & 1 \end{bmatrix}^\intercal \in \mathbb{R}^{c \times 1}$. Then we have
    \[
        F_p = \sum_{i = 1}^c (\bm{1}^\intercal x_i^\intercal x_i \bm{1}) \, \lambda_i^p \qquad (p \in \mathbb{N}_0).
    \]
\end{Cor}

\end{document}